\documentclass{article}

\usepackage{graphicx}
\usepackage{amsmath}
\usepackage{amsfonts}
\usepackage{amssymb}
\usepackage{subfigure}
\usepackage[amssymb]{SIunits}
\usepackage{color}

\newcommand{\Eq}[1]{(\ref{#1})}
\newcommand{\V}[1]{\ensuremath{\mathbf{#1}}}
\newcommand{\Vg}[1]{\boldsymbol{#1}}
\newcommand{\GV}[1]{\ensuremath{\mathbf{#1}}}                
\newcommand{\GM}[1]{\ensuremath{\mathbf{#1}}}           
\newcommand{\dtn}{\Delta t_n \,}

\title{Fast Solvers for Unsteady Thermal Fluid Structure Interaction}

\author{Philipp Birken$^{\times *}$, Tobias Gleim$^{\dag}$, Detlef Kuhl$^{\dag}$ and Andreas Meister$^{*}$}

\begin{document}
\maketitle

\begin{center}{\footnotesize\em $^{\times}$Centre for the Mathematical
    Sciences, Numerical Analysis, Lunds University, Box 118, 22100 Lund, Sweden\\
e-mail: philipp.birken@na.lu.se\\
$^{*}$Institute of Mathematics, University of Kassel, Heinrich-Plett-Str. 40, 34132 Kassel, Germany\\
$^{\dag}$ Institute of Mechanics and Dynamics, University of Kassel, M\"onchebergstr. 7, 34109 Kassel, Germany}
\end{center}

\begin{abstract}
We consider time dependent thermal fluid structure interaction. The respective models
are the compressible Navier-Stokes equations and the nonlinear heat
equation. A partitioned coupling approach via a
Dirichlet-Neumann method and a fixed point iteration is employed. As a refence solver a previously developed efficient
time adaptive higher order time integration scheme is used.

To improve upon this, we work on reducing the number of fixed point
coupling iterations. Thus, first widely used vector extrapolation methods for
convergence acceleration of the fixed point iteration are tested. In particular, Aitken relaxation, minimal polynomial
extrapolation (MPE) and reduced rank extrapolation (RRE) are
considered. Second, we explore the idea of extrapolation based on data
given from the time integration and derive such methods for SDIRK2. While the vector extrapolation
methods have no beneficial effects, the extrapolation methods allow to
reduce the number of fixed point iterations further by up to a factor
of two with linear extrapolation performing better than quadratic. 
\end{abstract}

{\it Keywords: Thermal Fluid Structure Interaction, Partitioned
  Coupling, Convergence Acceleration, Extrapolation}

\section{Introduction}

Thermal interaction between fluids and structures plays an important
role in many applications. Examples for this are cooling of gas-turbine blades, thermal anti-icing systems of airplanes \cite{buchli:10} or supersonic reentry of vehicles from space \cite{mehta:05,hinrad:06}. Another is quenching, an industrial heat treatment of metal workpieces. There, the desired material properties are achieved by rapid local cooling, which causes solid phase changes, allowing to create graded materials with precisely defined properties. 

Gas quenching recently received a lot of industrial and scientific interest \cite{wesast:07,hefiba:01}. In contrast to liquid quenching, this process has the advantage of minimal environmental impact because of non-toxic quenching media and clean products like air \cite{stshle:06}. To exploit the multiple advantages of gas quenching the application of computational fluid dynamics has proved essential \cite{banka:05,stshle:06,lior:04}. Thus, we consider the coupling of the compressible Navier-Stokes equations as a model for air, along a non-moving boundary with the nonlinear heat equation as a model for the temperature distribution in steel.

For the solution of the coupled problem, we prefer a partitioned approach
\cite{farhat:04}, where different codes for the sub-problems are reused
and the coupling is done by a master program which calls interface
functions of the other codes. This allows to use existing software for
each sub-problem, in contrast to a monolithic approach, where a new
code is tailored for the coupled equations. To satisfy the boundary
conditions at the interface, the subsolvers are iterated in a fixed
point procedure. 

Our goal here is to find a fast solver in this partitioned
setting. One approach would be to speed up the subsolvers and there is active research on that. See \cite{birkenhabil} for the current situation for fluid solvers. However, we want to
approach the problem from the point of view of a partitioned coupling method,
meaning that we use the subsolvers as they are. As a reference solver,
we use the time adaptive higher order time integration method suggested in \cite{biquhm:11}. Namely, the singly diagonally implicit Runge-Kutta  (SDIRK) method SDIRK2 is employed. 

To improve upon this, one idea is to define the tolerances
in the subsolver in a smart way and recently, progress has been made
for steady problems \cite{birken:14}. However, it is not immediately
clear how to transfer these results to the unsteady case. Thus, the
most promising way is to reduce the number of fixed point iterations,
on which we will focus in the present article. 

Various methods have been proposed to increase the convergence speed of the
fixed point iteration by decreasing the interface error between
subsequent steps, for example Relaxation \cite{letmou:01,kuewal:08},
Interface-GMRES \cite{mivabo:06}, ROM-coupling
\cite{viladv:06} and multigrid coupling \cite{zubobi:07}. Here,
we consider the most standard method, namely Aitken Relaxation and two variants
of polynomial vector extrapolation, namely MPE and RRE
\cite{sidi:12}. These have the merit of being purely algebraic and
very easy to implement. 

The second idea we follow is that of extrapolation based on
knowledge about the time integration scheme. This has been successfully used in other contexts
\cite{arnold:10,erbdue:12}, but has to our knowledge never been tried in Fluid Structure Interaction, where typically little attention is given to the time integration. Here, we use linear and quadratic
extrapolation of old values from the time history, designed specifically
for SDIRK2. 

The various methods are compared on the basis of numerical examples, namely the flow
past a flat plate, a basic test case for thermal fluid structure
interaction and an example from gas quenching \cite{wesast:07}.

\section{Governing Equations}

The basic setting we are in is that on a domain $\Omega_1 \subset
\mathbb{R}^d$ the physics is described by a fluid model, whereas on a
domain $\Omega_2 \subset \mathbb{R}^d$, a different model describing the structure is used. The two domains are almost disjoint in that they
are connected via an interface. The part of the interface where the
fluid and the structure are supposed to interact is called the
coupling interface 
$\Gamma \subset \partial \Omega_1 \cup \partial \Omega_2$. Note that $\Gamma$ might be a true subset of the intersection, because the structure could be insulated. At the interface $\Gamma$, coupling conditions are prescribed that model the interaction between fluid and structure. For the thermal coupling problem, these conditions are that temperature and the normal component of the heat flux are continuous across the interface. 

\subsection{Fluid Model}
We model the fluid using the Navier-Stokes equations, which are a second order system of conservation laws (mass, momentum, energy) modeling viscous compressible flow. We consider the two dimensional case, written in conservative variables density $\rho$, momentum ${\bf m}=\rho {\bf v}$ and energy per unit volume $\rho E$:

\vskip-.6cm
\begin{eqnarray}
\partial_t \rho + \nabla \cdot {\textbf m} & = & 0, \nonumber \\
\partial_t m_i + \sum_{j=1}^2 \partial_{x_j} (m_i v_j + p \delta_{ij}) & = & \frac{1}{Re}\sum_{j=1}^2 \partial_{x_j} S_{ij}, \qquad i=1,2, \nonumber \\
\partial_t (\rho E) + \nabla \cdot  (H {\textbf m}) & = & \frac{1}{Re}\sum_{j=1}^2 \partial_{x_j} \left ( \sum_{i=1}^2 S_{ij} v_i - \frac{1}{Pr}q_j \right).\nonumber
\end{eqnarray}

Here, $\textbf{S} = (S_{ij})_{i,j=1,2}$ represents the viscous shear stress tensor and
${\bf q} = (q_1,q_2)^T$ the heat flux. As the equation are dimensionless, the
Reynolds number $Re$ and the Prandtl number $Pr$ appear. The system is closed by the equation of state for the pressure $p = (\gamma -1)
\rho e$, the Sutherland law representing the correlation between temperature and viscosity as well as the Stokes hypothesis. Additionally, we prescribe appropriate boundary conditions at
the boundary of $\Omega_1$ except for $\Gamma$, where we have the coupling conditions. In the Dirichlet-Neumann coupling, a temperature value is enforced strongly at $\Gamma$. 

\subsection{Structure Model}
Regarding the structure model, we will consider heat conduction
only. Thus, we have the nonlinear
heat equation for the structure temperature $\Theta$
\vskip-.6cm
\begin{eqnarray}
  \rho({\bf x}) c_p(\Theta) \frac{d}{dt}\Theta(\V{x},t) = - \nabla \cdot \V{q}(\V{x},t),
  \label{eq:heatcond}
\end{eqnarray}
where 
\begin{equation*}
  \V{q}(\V{x},t) = - \lambda(\Theta) \nabla \Theta(\V{x},t)
\end{equation*}
denotes the heat flux vector. For steel, the specific heat capacity $c_p$ and heat conductivity $\lambda$ are temperature-dependent and highly nonlinear. Here, an empirical model for the steel 51CrV4 suggested in \cite{quhrss:11} is used. This model is characterized by the coefficient functions
\vskip-.6cm
\begin{eqnarray}
\lambda(\Theta)=40.1 + 0.05 \Theta - 0.0001\Theta^2 + 4.9 \cdot 10^{-8}\Theta^3
\end{eqnarray}
and
\vskip-.6cm
\begin{eqnarray}
c_p(\Theta)=-10\ln \left( \frac{e^{-c_{p1}(\Theta)/10}+e^{-c_{p2}(\Theta)/10}}{2}\right)
\end{eqnarray} 
with 
\vskip-.6cm
\begin{eqnarray}
c_{p1}(\Theta)=34.2e^{0.0026\Theta} + 421.15
\end{eqnarray} 
and 
\vskip-.6cm
\begin{eqnarray}
c_{p2}(\Theta)=956.5e^{-0.012(\Theta-900)} + 0.45\Theta.
\end{eqnarray} 

For the mass density one has $\rho =
\unit{7836}{\kilo\gram\per\cubic\meter}$. Finally, on the boundary, we have Neumann
conditions $\V{q}(\V{x},t) \cdot \V{n}(\V{x}) =
q_b(\V{x},t)$.

\section{Discretization}

\subsection{Discretization in space}
Following the partitioned coupling approach, we discretize the two models separately in space. For the fluid, we use a finite volume method, leading to 

\begin{eqnarray}
\label{eq:ODEfluid}
\frac{d}{dt} {\bf u} + {\bf h}({\bf u},\Vg{\Theta}) = \V{0},
\end{eqnarray}
where ${\bf h}({\bf u},\Vg{\Theta})$ represents the spatial discretization and its dependence on the temperatures in the fluid. In particular, the DLR TAU-Code is employed \cite{gerfeg:97}, which is a cell-vertex-type finite volume method with AUSMDV as flux function and a linear reconstruction to increase the order of accuracy. 

Regarding structural mechanics, the use of finite element methods is ubiquitious. Therefore, we will also follow that approach here and use quadratic finite elements \cite{zienkiewicztaylor00a}, leading to the nonlinear equation for all unknowns on $\Omega_2$

\begin{eqnarray}
\label{eq:ODEheat}
\GM{M}({\bf \Theta}) \frac{d}{dt}{\bf \Theta} + \GM{K}({\bf \Theta}) {\bf \Theta}
  = {\bf q}_b({\bf u}).
\end{eqnarray}
Here, $\GM{M}$ is the heat capacity and $\GM{K}$ the heat conductivity
matrix. The vector ${\bf \Theta}$ consists of all discrete temperature
unknowns and ${\bf q}_b$ is the heat flux vector on the
surface. In this case it is the prescribed Neumann heat flux vector of
the fluid.

\subsection{Coupled time integration}

For the time integration, a time adaptive SDIRK2 method is implemented in a partitioned way, as suggested in \cite{biquhm:11}. If the fluid and the solid solver are able to carry out time steps of implicit Euler type, the master program of the FSI procedure can be extended to SDIRK methods very easily, since the master program just has to call the backward Euler routines with specific time step sizes and starting vectors. This method is very efficient and will be used as the base method in its time adaptive variant, which is much more efficient than more commonly used fixed time step size schemes. 

To obtain time adaptivity, embedded methods are used. Thereby, the local error is estimated by the solvers
separately, which then report the estimates back to the master program. Based on this, the new time step is chosen \cite{biquhm:10}. To this end, all stage derivatives are stored by the subsolvers. Furthermore, if the possibility of rejected time steps is taken into account, the current solution pair $({\bf u}, {\bf \Theta})$ has to be stored as well. 

To comply with the conditions that the discrete temperature and heat flux are
continuous at the interface $\Gamma$, a Dirichlet-Neumann
coupling is used. Thus, the boundary conditions for the two solvers
are chosen such that we prescribe Neumann data for one solver and
Dirichlet data for the other. Following the analysis of Giles
\cite{giles:97}, temperature is prescribed for the equation with
smaller heat conductivity, here the fluid, and heat flux is given on
$\Gamma$ for the structure. Choosing these conditions the other way
around leads to an unstable scheme. 

In the following it is assumed that at time $t_n$, the step size $\dtn$ is prescribed. Applying a DIRK method to equation \Eq{eq:ODEfluid}-\Eq{eq:ODEheat} results in the coupled system of equations to be solved at Runge-Kutta stage $i$:
\vskip-.6cm
\begin{eqnarray}
  \label{eq:Fnoneq}
  \mathbf{F}({\bf u}_i,\Vg{\Theta}_i) :=
  {\bf u}_i - {\bf s}_i^{{\bf u}} - \dtn a_{ii}
  {\bf h}({\bf u}_i,\Vg{\Theta}_i) = \V{0},
\end{eqnarray}
\vskip-.6cm
\begin{eqnarray}\label{eq:Tlineq}
  \V{T}({\bf u}_i,\Vg{\Theta}_i) :=
    [\GM{M} - \dtn a_{ii} \GM{K}] \Vg{\Theta}_i -
    \GM{M} \V{s}_i^{\Theta} - \dtn a_{ii}\V{q}_b({\bf u}_i) =\V{0}.
  \end{eqnarray}
Here, $a_{ii}=1-\sqrt{2}/2$ is a coefficient of the time integration method and ${\bf s}_i^{{\bf u}}$ and $\V{s}_i^{\Theta}$ are given vectors, called starting vectors, computed inside the DIRK scheme.
The dependence of the fluid equations ${\bf h}({\bf u}_i,\Vg{\Theta}_i)$ on the temperature $\Vg{\Theta}_i$ results from the nodal temperatures of the structure at the interface. This subset is written as $\Vg{\Theta}_i^{\Gamma}$. Accordingly, the structure equations depend only on the heat flux of the fluid at the coupling interface.

\section{Fixed Point Iteration and Improvements}

\subsection{Basic fixed point iteration}
To solve the coupled system of nonlinear equations
\Eq{eq:Fnoneq}-\Eq{eq:Tlineq}, a strong
coupling approach is employed. Thus, a fixed point iteration
is iterated until a convergence criterion is satisfied. In particular,
we use a a nonlinear Gau\ss-Seidel process:
\begin{align}
  \V{F}({\bf u}_i^{(\nu+1)},\Vg{\Theta}_i^{(\nu)})
  &=
  \GV{0} 
  \quad \leadsto \quad {\bf u}_i^{(\nu+1)}
  \label{eq:GSF}\\
  \V{T}({\bf u}_i^{(\nu+1)},\Vg{\Theta}_i^{(\nu+1)})
  &=
  \GV{0} 
  \quad \leadsto \quad \Vg{\Theta}_i^{(\nu+1)}, \quad \nu=0,1,...
  \label{eq:GST}
\end{align}
Each inner iteration is thereby done locally by
the structure or the fluid solver. More specific, a Newton method is
used in the structure and a
FAS multigrid method is employed in the fluid. 

In the base method, the starting
values of the iteration are given by ${\bf u}_i^{(0)}=\V{s}_i^{{\bf u}}$ and
$\Vg{\Theta}_i^{(0)}=\V{s}_i^{\Theta}$. The termination criterion is
formulated by the relative update of the nodal temperatures at the interface of the solid structure and we stop once we are below the tolerance in the time integration scheme divided by five
\vskip-.6cm
\begin{eqnarray}
  \label{eq:FP}
  \|\Vg{\Theta}_i^{\Gamma \; (\nu+1)} - \Vg{\Theta}_i^{\Gamma \;
    (\nu)}\|
  \leq \mathit{TOL}/5\|\Vg{\Theta}_i^{\Gamma \; (0)}\|.
\end{eqnarray}
The vector
\vskip-.6cm
\begin{eqnarray}
  \label{eq:residual}
{\bf r}^{(\nu +1)}:=\Vg{\Theta}_i^{\Gamma \; (\nu+1)} - \Vg{\Theta}_i^{\Gamma \;
    (\nu)}
\end{eqnarray}
is often referred to as the interface residual. 

We will now consider different techniques to improve upon this base iteration,
namely using vector extrapolation inside the fixed point iteration and
then extrapolation inside the time integration schemes, to obtain
better initial values. 

\subsection{Vector Extrapolation}
To improve the convergence speed of the fixed point iteration, different vector extrapolation techniques have been suggested. These are typically classic techniques, where a set of $k$ vectors of a convergent vector sequence
is extrapolated to obtain a faster converging sequence. We are now
going to describe three techniques that we will investigate in this
framework.

\subsubsection{Aitken Relaxation}

Relaxation means that after the fixed point iterate is computed, a relaxation step is added:

\vskip-.6cm
\begin{eqnarray}
  \label{eq:relax_a}
  \tilde{\Vg{\Theta}}_i^{\Gamma \; (\nu+1)} = \omega_{\nu +1}
  \Vg{\Theta}_i^{\Gamma \; (\nu+1)} + (1-\omega_{\nu +1}) \Vg{\Theta}_i^{\Gamma \; (\nu)}.
  \end{eqnarray}
Several strategies exist to compute the relaxation parameter
$\omega$.

The idea of Aitken's method is to enhance the current solution
$\Vg{\Theta}_i^{\Gamma \; (\nu+1)}$ using two previous iteration pairs
$(\Vg{\Theta}_i^{\Gamma \; (\nu+2)}, \Vg{\Theta}_i^{\Gamma \;
  (\nu+1)})$ and $(\Vg{\Theta}_i^{\Gamma \; (\nu+1)},
\Vg{\Theta}_i^{\Gamma \; (\nu)})$ obtained from the Gau\ss -Seidel-step \Eq{eq:GSF}-\eqref{eq:GST}. An improvement in the scalar case is given by the secant method
\vskip-.6cm
\begin{eqnarray}
  \label{eq:secant}
    \tilde{{\Theta}}_i^{\Gamma \; (\nu+1)} = \dfrac{{\Theta}_i^{\Gamma \; (\nu-1)} {\Theta}_i^{\Gamma \; (\nu+1)} - {\Theta}_i^{\Gamma \; (\nu)} {\Theta}_i^{\Gamma \; (\nu)}}{{\Theta}_i^{\Gamma \; (\nu-1)} - {\Theta}_i^{\Gamma \; (\nu)} - {\Theta}_i^{\Gamma \; (\nu)} + {\Theta}_i^{\Gamma \; (\nu+1)}}.
\end{eqnarray}
The relaxation factor in equation \Eq{eq:relax_a} for the secant method \Eq{eq:secant} is then
\vskip-.6cm
\begin{eqnarray}
  \label{eq:relax_b}
    \omega_{\nu+1} = \dfrac{{\Theta}_i^{\Gamma \; (\nu-1)} -  {\Theta}_i^{\Gamma \; (\nu)}}{{\Theta}_i^{\Gamma \; (\nu-1)} - {\Theta}_i^{\Gamma \; (\nu)} -      {\Theta}_i^{\Gamma \; (\nu)} + {\Theta}_i^{\Gamma \; (\nu+1)}}.
\end{eqnarray}
As customary, we use an added recursion on $\omega_i$ in which we use the old relaxation factor $\omega_{\nu}$:
\vskip-.6cm
\begin{eqnarray}
  \label{eq:omega_recursion}
    \omega_{\nu+1} = -\omega_{\nu}\ \dfrac{{r}^{\Gamma \;  (\nu)}}{{ r}^{\Gamma \;  (\nu+1)} - {r}^{\Gamma \;  (\nu)}}. 
\end{eqnarray}
In the vector case the division by the residual ${\bf r}^{\Gamma \;  (\nu+1)} - {\bf r}^{\Gamma \;  (\nu)}$ is not possible. Therefore, we multiply the nominator and the numerator formally by $({\bf r}^{\Gamma \;  (\nu+1)} - {\bf r}^{\Gamma \;  (\nu)})^T$ to obtain
\vskip-.6cm
\begin{eqnarray}
  \label{eq:relax_c}
    \omega_{\nu+1} = -\omega_{\nu}\frac{({\bf r}^{\Gamma \;  (\nu)})^T({\bf r}^{\Gamma \;  (\nu +1)}-{\bf r}^{\Gamma \;  (\nu)})}{\|{\bf r}^{\Gamma \;  (\nu +1)}-{\bf r}^{\Gamma \;  (\nu)}\|^2}.
\end{eqnarray}
Two previous steps are required to calculate the relaxation parameter. For the first fixpoint iteration, the relaxation parameter $\omega_1$ must be prescribed. We choose $\omega_1= 0.8$, which was reported by other authors to work well. 

\subsubsection{Polynomial Vector Extrapolation}
Another idea we will follow here are Minimal Polynomial Extrapolation (MPE) and Reduced Rank Extrapolation (RRE) \cite{sidi:12}. Here, the new approximation is given as a linear combination of existing iterates with coefficients $\gamma_{\nu}$ to be determined:
\vskip-.6cm
\begin{eqnarray}
  \label{eq:vec_approx}
    \tilde{\Vg{\Theta}}_i^{\Gamma \; (\nu+1)} =
    \sum\limits_{j=0}^{\nu +1} \gamma_j \Vg{\Theta}_i^{\Gamma \; (j)}.
\end{eqnarray}
For MPE, the coefficients are defined via
\vskip-.6cm
\begin{eqnarray}
\gamma_j=\frac{c_j}{\sum_{i=0}^{\nu+1}c_i}, \quad j=0,...,\nu+1,
\end{eqnarray}
where the coefficients $c_j$ are the solution of the problem
\vskip-.6cm
\begin{eqnarray}
\min_{c_j}\Big \|\sum_{j=0}^{\nu+1}c_j{\bf r}_j+{\bf r}_{\nu+1}\Big \|_2.
\end{eqnarray}

For RRE, the coefficients are defined as the solution of the constrained least squares problem
\vskip-.6cm
\begin{eqnarray}
\min_{\gamma_j}\Big \|\sum_{j=0}^{\nu+1}\gamma_j{\bf r}_j\Big \|_2 \quad {\rm subject}\,{\rm to} \quad \sum_{j=0}^{\nu+1} \gamma_j=1.
\end{eqnarray}

These problems are then solved using a QR decomposition.

\subsection{Extrapolation from time integration}

To find good starting values for iterative processes in implicit time integration schemes, it is common to use extrapolation based on knowledge about the trajectory of the solution of the initial value problem \cite{erbdue:12,olssoe:98}. In the spirit of partitioned solvers, we here suggest to use extrapolation of the interface temperatures only. Of course, this strategy could be used as well by the subsolvers, but we will not consider this here. 

At the first stage, we have the old time step size $\Delta t_{n-1}$ with value $\Theta_{n-1}$ and the current time step size $\Delta t_n$ with value $\Theta_n$. We are looking for the value $\Theta_1$ at the next stage time $t_n + c_1\Delta t_n$. Linear extrapolation results in
\begin{equation}
\Theta_1=\Theta_n + c_1\Delta t_n(\Theta_n-\Theta_{n-1})/\Delta t_{n-1}=(1+\frac{c_1\Delta t_n}{\Delta t_{n-1}})\Theta_n-\frac{c_1\Delta t_n}{\Delta t_{n-1}}\Theta_{n-1}.
\end{equation}
Regarding quadratic extrapolation, it is reasonable to choose $t_n$, $t_{n-1}$ and the intermediate temperature vector $\Theta_{n-1/2}$ from the previous stage $t_{n-1} + c_1\Delta t_{n-1}$. This results in
\begin{eqnarray}
\Theta_1=&\Theta_{n-1}\frac{(c_1\Delta t_n+(1-c_1)\Delta t_{n-1})c_1\Delta t_n}{c_1\Delta t_{n-1}^2} - \Theta_{n-1/2}\frac{(c_1\Delta t_n+\Delta t_{n-1})c_1\Delta t_n}{c_1\Delta t_{n-1}^2(1-c_1)} \nonumber \\
&+ \Theta_n\frac{(c_1\Delta t_n + \Delta t_{n-1})(c_1\Delta t_n + (1-c_1)\Delta t_{n-1})}{(1-c_1)\Delta t_{n-1}^2}.
\end{eqnarray}

When applying this idea at the second stage (or at later stages in a scheme with more than two),  it is better to use values from the current time interval. Thus, we linearly extrapolate $\Theta_n$ at $t_n$ and $\Theta_1$ at $t_n + c_1\Delta t$ to obtain
\begin{equation}
\label{eq:lin_extrap}
\Theta_{n+1}=\Theta_n + \Delta t_n(\Theta_1-\Theta_n)/(c_1\Delta t_n)=(1-\frac{1}{c_1})\Theta_n+\frac{1}{c_1}\Theta_1.
\end{equation}
This results in
\begin{eqnarray}
\label{eq:quad_extrap}
\Theta_{n+1}=&\Theta_{n-1}\frac{\Delta t_n^2(1-c_1)}{\Delta t_{n-1}(\Delta t_{n-1}+c_1\Delta t_n)} - \Theta_n\frac{(\Delta t_{n-1}+\Delta t_n)(1-c_1)\Delta t_n}{\Delta t_{n-1}c_1\Delta t_n} \nonumber \\
& +\Theta_1\frac{(\Delta t_{n-1}+\Delta t_n)\Delta t_n}{(c_1\Delta t_n+\Delta t_{n-1})c_1\Delta t_n}.
\end{eqnarray}

\section{Numerical Results}

\subsection{Flow over a plate}
\begin{figure}[ht]
  \includegraphics[width=\linewidth]{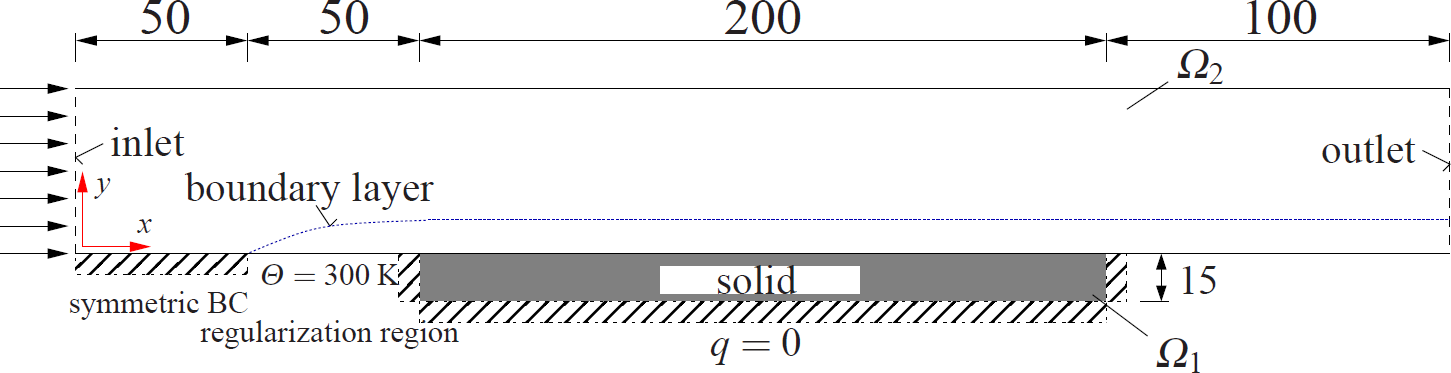}
 \caption{Test case for the coupling method}
 \label{fig:testcase}
\end{figure}

As a first test case, the cooling of a flat plate resembling a simple work piece is considered. The work piece is initially at a much higher temperature than the fluid and then cooled by a constant air stream, that is assumed to be laminar, see figure \ref{fig:testcase}.

The inlet is given on the left, where air enters the domain with an
initial velocity of ${\rm Ma}_{\infty} = 0.8$ in horizontal direction and a temperature of \unit{273}{\kelvin}. Then, there are two succeeding regularization regions of \unit{50}{\milli\meter} to obtain an unperturbed boundary layer. In the first region, $0 \leq x \leq 50$, symmetry boundary conditions, $v_y=0$, $q=0$, are applied. In the second region, $50 \leq x \leq 100$, a constant wall temperature of \unit{300}{\kelvin} is specified. Within this region the velocity boundary layer fully develops. The third part is the solid (work piece) of length \unit{200}{\milli\meter}, which exchanges heat with the fluid, but is assumed insulated otherwise, thus $q_b=0$. Therefore, Neumann boundary conditions are applied throughout. Finally, the fluid domain is closed by a second regularization region of \unit{100}{\milli\meter} with symmetry boundary conditions and the outlet. 

Regarding the initial conditions in the structure, a constant temperature of \unit{900}{\kelvin} at $t=\unit{0}{\second}$ is chosen throughout. To specify reasonable initial conditions within the fluid, a steady state solution of the fluid with a constant wall temperature $\Theta=\unit{900}{\kelvin}$ is computed. 

The grid
\begin{figure}[h!]
  \centering
  \subfigure[Entire mesh]{
    \includegraphics[width=0.7
    \linewidth]{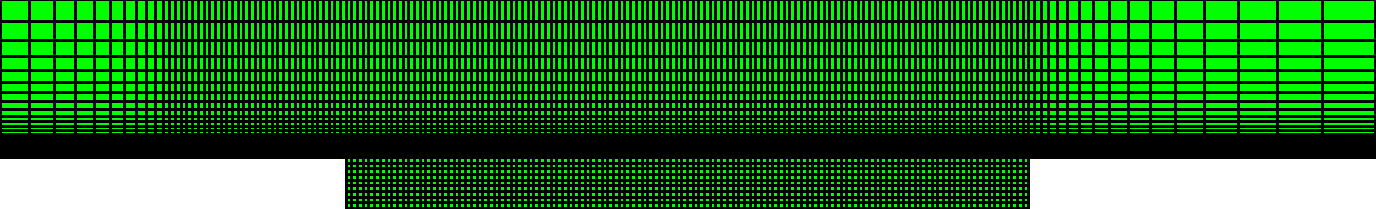}}
\hfill
  \subfigure[Mesh zoom]{
    \includegraphics[width=0.25 \linewidth]{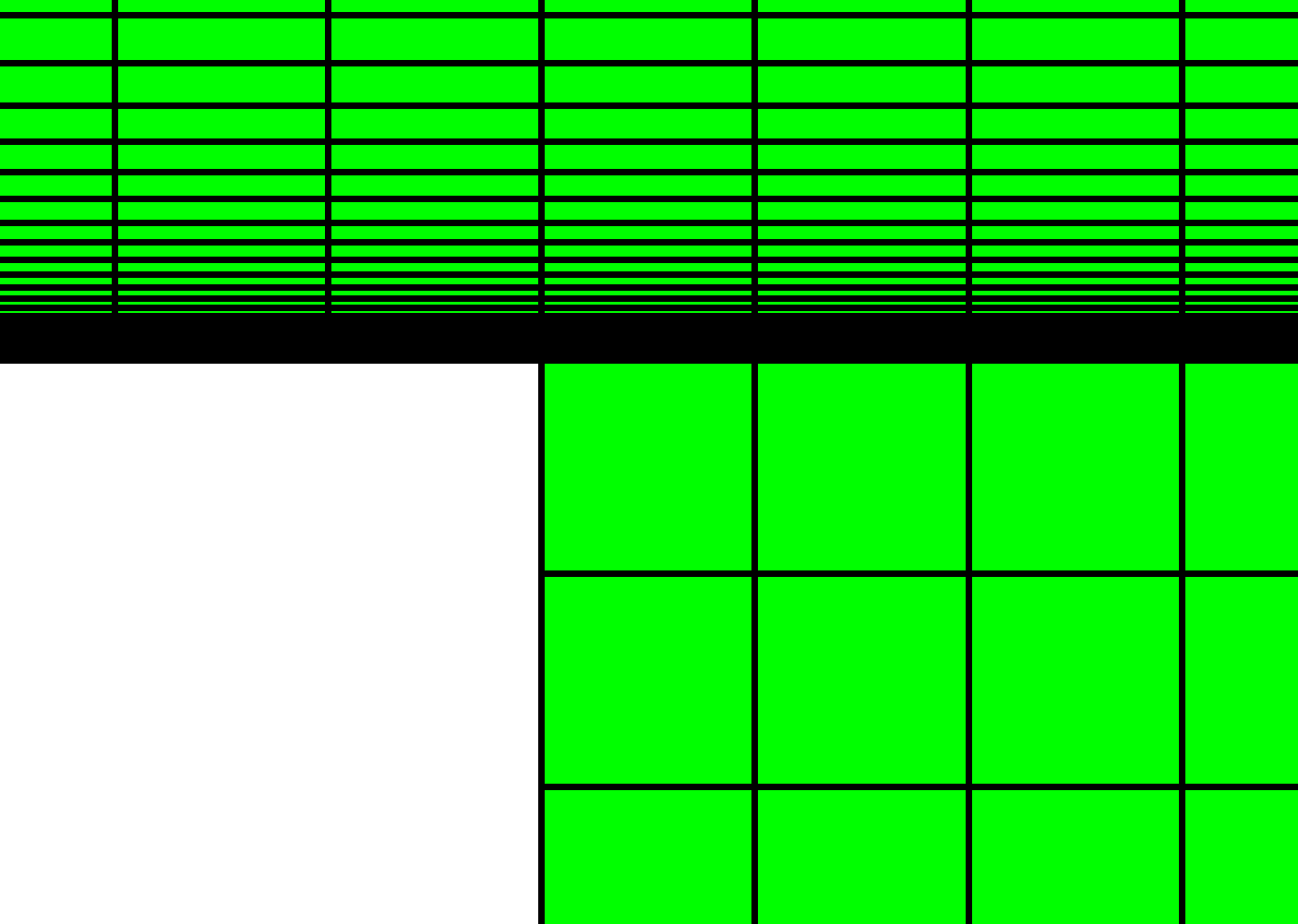}}
  \caption{Full grid (left) and zoom into coupling region (right)}
  \label{fig:grid}
\end{figure}
is chosen cartesian and equidistant in the structural part, where in the fluid region the thinnest cells are on the boundary and then become coarser in $y$-direction (see figure \ref{fig:grid}). To avoid additional difficulties from interpolation, the points of the primary fluid grid, where the heat flux is located in the fluid solver, and the nodes of the structural grid are chosen to match on the interface $\Gamma$. 

\newcommand{\lokal}{\scriptsize}
\newcommand{\lokalc}{\footnotesize}
\newcommand{\lokalb}{\tiny}
\begin{figure}[t]
\centering
 \subfigure[1 time step with a time step size of \unit{1}{s}]{
 \begin{picture}(160,140)
 \put(2,0){
\put(20,15){\includegraphics[width=0.42\textwidth]{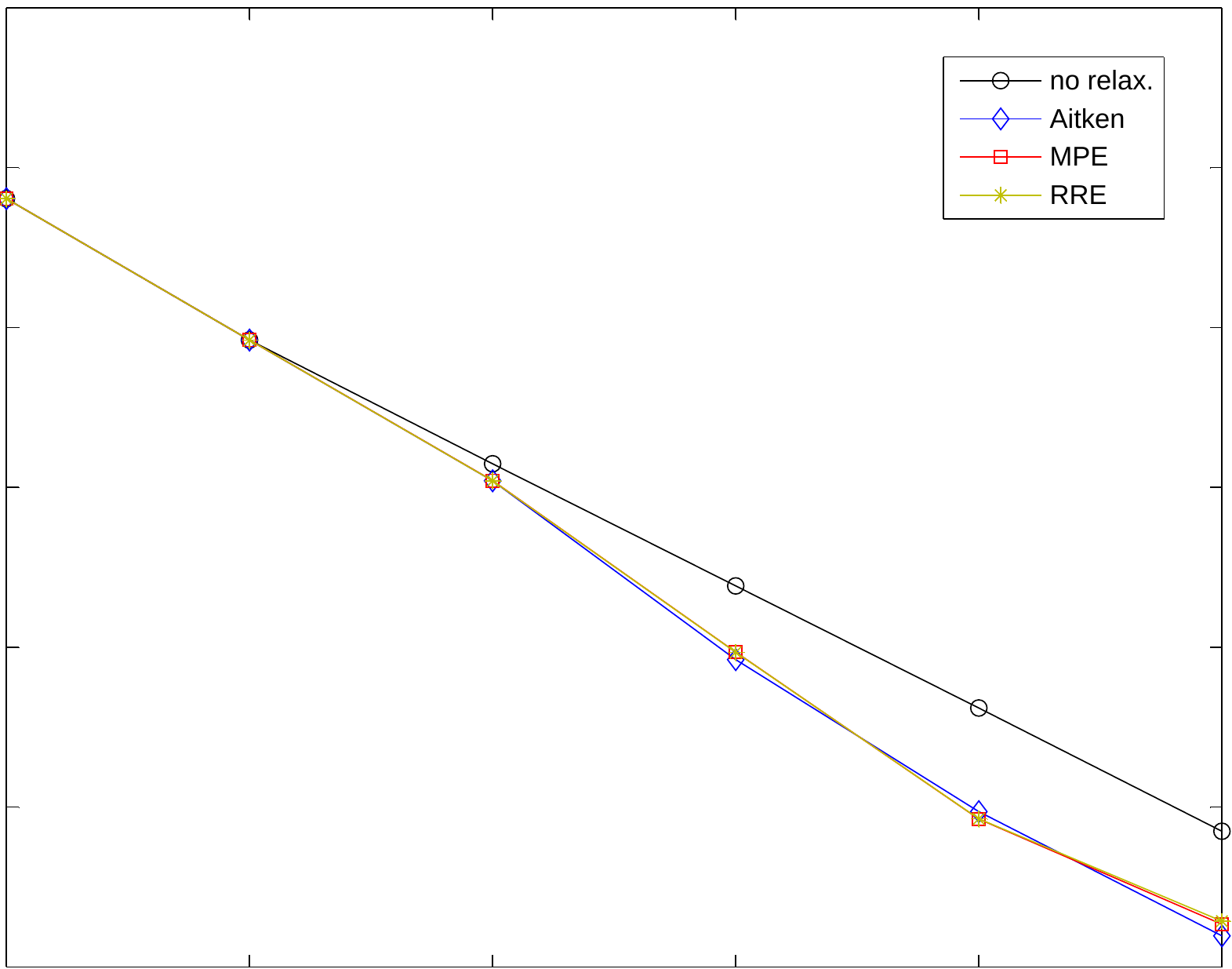}}
\put(120,100){\textcolor{white}{\rule[0mm]{1.4cm}{0.9cm}}}
\put(0,-2){
\put(100,120){\includegraphics[width=5mm]{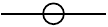}}
\put(100,110){\includegraphics[width=5mm]{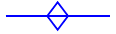}}
\put(100,100){\includegraphics[width=5mm]{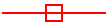}}
\put(100,90){\includegraphics[width=5mm]{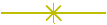}}
         \put(115,123){\lokal{\makebox(0,0)[l]{no relax.}}} 
                  \put(115,112){\lokal{\makebox(0,0)[l]{Aitken}}} 
                           \put(115,102){\lokal{\makebox(0,0)[l]{MPE}}} 
                                    \put(115,92){\lokal{\makebox(0,0)[l]{RRE}}} 
                                    }
 \put(21,10){\lokal{\makebox(0,0)[c]{$1$}}}
 \put(50,10){\lokal{\makebox(0,0)[c]{$2$}}}
  \put(78,10){\lokal{\makebox(0,0)[c]{$3$}}}
   \put(107,10){\lokal{\makebox(0,0)[c]{$4$}}}
    \put(135,10){\lokal{\makebox(0,0)[c]{$5$}}}
     \put(163,10){\lokal{\makebox(0,0)[c]{$6$}}} 
          \put(90,3){\lokal{\makebox(0,0)[c]{Iteration}}} 
          \put(10,0){
           \put(9,18){\lokalb{\makebox(0,0)[r]{$10^{-12}$}}}
           \put(9,36){\lokalb{\makebox(0,0)[r]{$10^{-10}$}}}
               \put(9,54){\lokalb{\makebox(0,0)[r]{$10^{-8}$}}}
                         \put(-11,75){\lokalc{\makebox(0,0)[r]{\rotatebox{90}{2-norm of interface residual}}}} 
                   \put(9,72){\lokalb{\makebox(0,0)[r]{$10^{-6}$}}}
                       \put(9,91){\lokalb{\makebox(0,0)[r]{$10^{-4}$}}}   
                           \put(9,110){\lokalb{\makebox(0,0)[r]{$10^{-2}$}}}
                               \put(9,128){\lokalb{\makebox(0,0)[r]{$10^{0}$}}}}     
                               }
\end{picture}}
\hfill
 \subfigure[1 time step with a time step size of \unit{5}{s}]{
  \begin{picture}(160,140)
  \put(0,0){
\put(20,15){\includegraphics[width=0.42\textwidth]{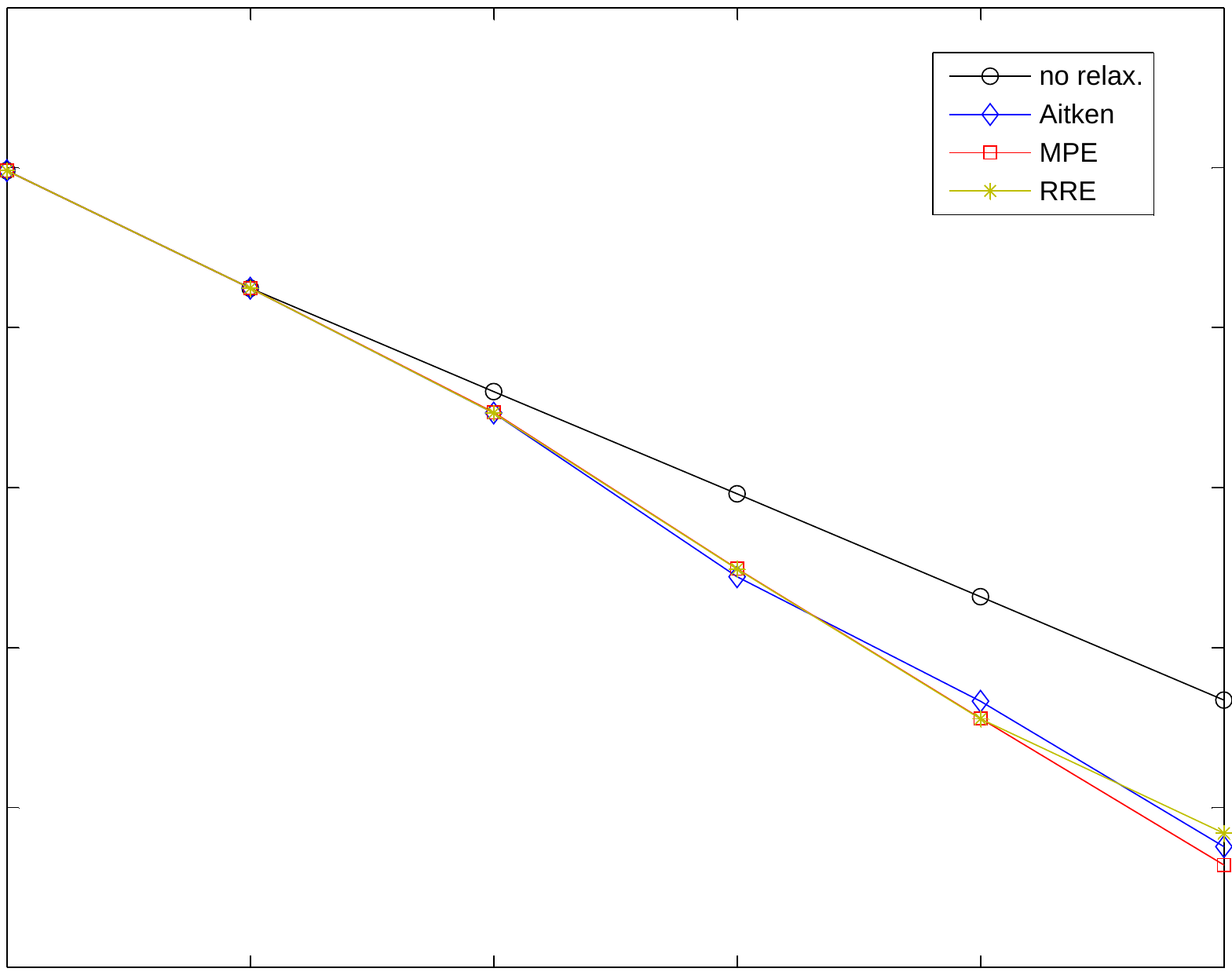}}
\put(120,100){\textcolor{white}{\rule[0mm]{1.4cm}{0.9cm}}}
\put(0,-2){
\put(100,120){\includegraphics[width=5mm]{legende_kreis.pdf}}
\put(100,110){\includegraphics[width=5mm]{legende_raute.pdf}}
\put(100,100){\includegraphics[width=5mm]{legende_rule.pdf}}
\put(100,90){\includegraphics[width=5mm]{legende_stern.pdf}}
         \put(115,123){\lokal{\makebox(0,0)[l]{no relax.}}} 
                  \put(115,112){\lokal{\makebox(0,0)[l]{Aitken}}} 
                           \put(115,102){\lokal{\makebox(0,0)[l]{MPE}}} 
                                    \put(115,92){\lokal{\makebox(0,0)[l]{RRE}}} 
                                    }
 \put(21,10){\lokal{\makebox(0,0)[c]{$1$}}}
 \put(50,10){\lokal{\makebox(0,0)[c]{$2$}}}
  \put(78,10){\lokal{\makebox(0,0)[c]{$3$}}}
   \put(107,10){\lokal{\makebox(0,0)[c]{$4$}}}
    \put(135,10){\lokal{\makebox(0,0)[c]{$5$}}}
     \put(163,10){\lokal{\makebox(0,0)[c]{$6$}}} 
          \put(90,3){\lokal{\makebox(0,0)[c]{Iteration}}} 
          \put(10,0){
           \put(9,18){\lokalb{\makebox(0,0)[r]{$10^{-12}$}}}
           \put(9,36){\lokalb{\makebox(0,0)[r]{$10^{-10}$}}}
               \put(9,54){\lokalb{\makebox(0,0)[r]{$10^{-8}$}}}
                         \put(-11,75){\lokalc{\makebox(0,0)[r]{\rotatebox{90}{2-norm of interface residual}}}} 
                   \put(9,72){\lokalb{\makebox(0,0)[r]{$10^{-6}$}}}
                       \put(9,91){\lokalb{\makebox(0,0)[r]{$10^{-4}$}}}   
                           \put(9,110){\lokalb{\makebox(0,0)[r]{$10^{-2}$}}}
                               \put(9,128){\lokalb{\makebox(0,0)[r]{$10^{0}$}}}}
                               }
\end{picture}
}
\caption{\label{fig:onesystem}Comparison of the relaxation methods for 1 time step with different timestep sizes}
\end{figure}
To compare the effect of the different vector extrapolation strategies, we
consider the fixed point equation within the first stage of the first time step in the test problem with a time step size of $\Delta t=1$s and $\Delta t=5$s. In figure \ref{fig:onesystem}, we can see how the interface residual decreases with the fixed point iterations. During the first two steps all methods have the same residual norm, since all methods need at least two iterations to start. For this example, the vector extrapolation methods outperform the standard scheme for tolerances below $10^{-5}$.

\begin{table}[h!]
\caption{\label{tab:plate-fixed-v-adaptive}Total number of iterations for 100 secs of
  real time without any extrapolation. Fixed timestepsizes versus adaptive steering.}
\centering
\begin{tabular}{c|cc|c}
\hline
$TOL$ & \multicolumn{2}{c|}{Fixed time step size} & Time adapt., $\Delta t_0=0.5s$\\
\hline
$10^{-2}$ &$\Delta t=5s$ & 64 & 31 \\
$10^{-3}$ & $\Delta t=5s$ & 82 & 39 \\
$10^{-4}$ &$\Delta t=0.5s$ & 802  & 106 \\
\hline
\end{tabular}
\label{tab:plate-fixed-timestep}
\end{table}

We now compare the different schemes for a whole simulation of $100$
seconds real time. If not mentioned otherwise, the initial time step
size is $\Delta t=0.5s$. To first give an impression on the effect of
the time adaptive method, we look at fixed time step versus adaptive
computations in tabular \ref{tab:plate-fixed-timestep}. Thus, the
different tolerances for the time adaptive case lead to different time
step sizes and tolerances for the nonlinear system over the course of
the algorithm, whereas in the fixed time step size, they steer only
how accurate the nonlinear systems are solved. For the fixed time step
case, we chose $\Delta t=0.5s$ and $\Delta t=5s$, which roughly
corresponds to an error of $10^{-2}$ and $10^{-3}$, respectively
$10^{-4}$. Thus, computations in one line of tabular
\ref{tab:plate-fixed-timestep} correspond to similar errors. As can be
seen, the time adaptive method is in the worst case a factor two
faster and in the best case a factor of eight. Thus the time adaptive
computation serves from now on as the base method for the construction
of a fast solver.



\begin{table}[h!]
\caption{\label{tab:100secs}Total number of iterations for different tolerances and different vector extrapolation strategies for 100 secs of real time for the flat plate test case.}
\centering
\begin{tabular}{*{6}{c}}
\hline
$TOL$ & No relax. & Aitken & MPE & RRE \\
\hline
$10^{-2}$ &31 &32 &31 &31\\
$10^{-3}$ &39 &38 &39 &39\\
$10^{-4}$ & 106 & 103 & 106 & 106\\
$10^{-5}$ &857   &736  &857 &857\\
\hline
\end{tabular}
\end{table}
The next computations demonstrate the effect of vector extrapolation. With increasing time the time adaptive algorithm chooses larger time step sizes. The base method needs more fixed point iterations in the end of the time interval, while the other methods have remained roughly constant. The total number of fixed point iterations is shown in tabular \ref{tab:100secs}. As we can see, only Aitken relaxation has an advantage over the base method and that only for a tolerance of $10^{-5}$. For larger tolerances, all the methods need roughly the same number of iterations, which is also confirmed in Figure \ref{fig:onesystem}, where all methods overlap for $\|{\bf r}\|_2\leq10^{-4}$.

Essentially, the interplay between the fixed point iteration and the time adaptive scheme results in only two fixed point iterations being necessary per time step (compare equation \eqref{eq:FP}). Thus, the vector extrapolation methods have no effect. 

\begin{table}[h!]
\caption{\label{tab:extrap-plate}Total number of iterations for 100 secs of
  real time with extrapolation}
\centering
\begin{tabular}{*{7}{c}}
\hline
$TOL$ & none & lin. & quad. \\
\hline
$10^{-2}$ & 31 &19 &25\\
$10^{-3}$ & 39 &31 &32\\
$10^{-4}$ & 106 &73 &77\\
$10^{-5}$ & 857 &415 &414\\
\hline
\end{tabular}
\end{table}

Finally, we consider extrapolation based on the time integration
scheme. In table \ref{tab:extrap-plate}, the total number of
iterations for 100 seconds of real time is shown. As can be seen,
linear extrapolation speeds up the computations between $20\%$ and $50
\%$. Quadratic extrapolation leads to a speedup between $15\%$ and
$50\%$ being overall less efficient than the linear extrapolation procedure. Overall, we are thus able to simulate 100 seconds of real time for this problem for an engineering tolerance using only 19 calls to fluid and the structure solver each. 

\begin{figure}[h!]
\centering
 \begin{picture}(160,140)
\put(20,15){\includegraphics[width=0.42\textwidth]{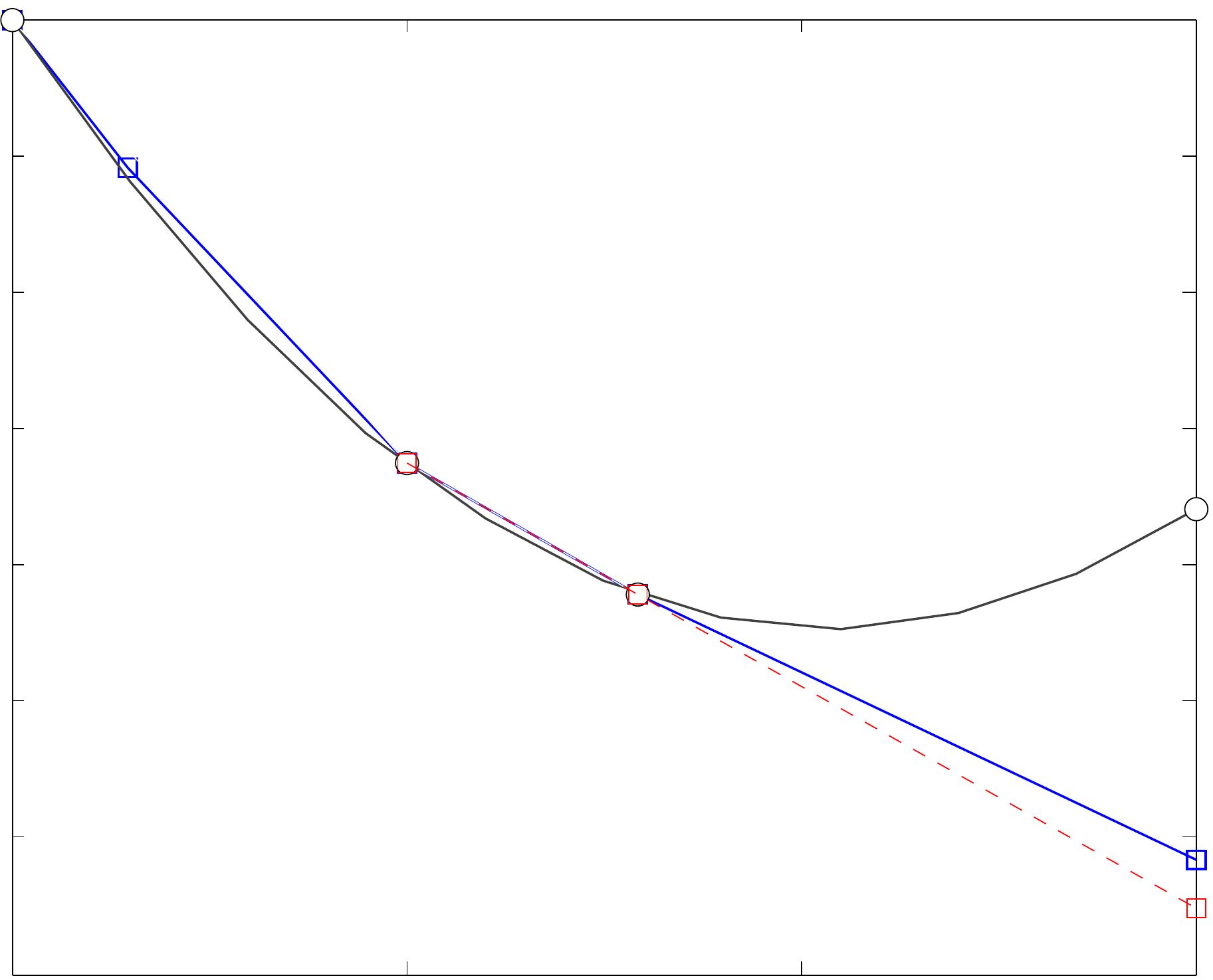}}
\put(120,100){\textcolor{white}{\rule[0mm]{1.4cm}{0.9cm}}}
\put(0,-2){
\put(100,120){\includegraphics[width=5mm]{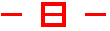}}
\put(100,110){\includegraphics[width=5mm]{legende_kreis.pdf}}
\put(100,100){\includegraphics[width=5mm]{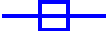}}
         \put(115,123){\lokal{\makebox(0,0)[l]{linear extrap.}}} 
                  \put(115,112){\lokal{\makebox(0,0)[l]{quad. extrap.}}} 
                           \put(115,102){\lokal{\makebox(0,0)[l]{final solns.}}} 
                                    }
 \put(21,10){\lokal{\makebox(0,0)[c]{$0$}}}
  \put(69,10){\lokal{\makebox(0,0)[c]{$0.5$}}}
   \put(116,10){\lokal{\makebox(0,0)[c]{$1$}}}
     \put(163,10){\lokal{\makebox(0,0)[c]{$1.5$}}} 
          \put(90,3){\lokal{\makebox(0,0)[c]{$t$ [s]}}} 
          \put(10,0){
           \put(9,18){\lokal{\makebox(0,0)[r]{$865$}}}
           \put(9,32){\lokal{\makebox(0,0)[r]{$870$}}}
               \put(9,50){\lokal{\makebox(0,0)[r]{$875$}}}
                   \put(9,66){\lokal{\makebox(0,0)[r]{$880$}}}
                   \put(9,82){\lokal{\makebox(0,0)[r]{$\Theta$ [K]}}}
                       \put(9,99){\lokal{\makebox(0,0)[r]{$890$}}}   
                           \put(9,115){\lokal{\makebox(0,0)[r]{$895$}}}
                               \put(9,130){\lokal{\makebox(0,0)[r]{$900$}}}}
\end{picture}

\caption{\label{fig:compare_extrapolation}Comparison of the linear and quadratic extrapolation methods for the time step $t = 1.5$s.}
\end{figure}

To understand this more precisely, we considered the second stage of the second time step in an
adaptive computation. We thus have finished the first time step with
$\Delta t_0=0.5s$ and the second time step gets doubled, leading to
$\Delta t_1=1s$. This is depicted in Figure
\ref{fig:compare_extrapolation}. 
To obtain a temperature for the new time $t_{n+1}$ the linear
extrapolation method \eqref{eq:lin_extrap} uses the values of the current time $t_n$ and of
the first Runge-Kutta Step at $t_1+\Delta t_1c_1$. As can be seen,
this predicts the new time step very well. In contrast, the quadratic
extrapolation \eqref{eq:quad_extrap} uses for the new time step the
solution from the previous time $t_0$ the current time $t_1$ and from
the first Runge Kutta stage. Since the exact solution has a more
linear behavior in the time step, the quadratic extrapolation provides
no advantage, in particular since it slopes upward after some
point. 

\subsection{Cooling of a flanged shaft}
\begin{figure}[h!]
\centering
  \includegraphics[width=60mm]{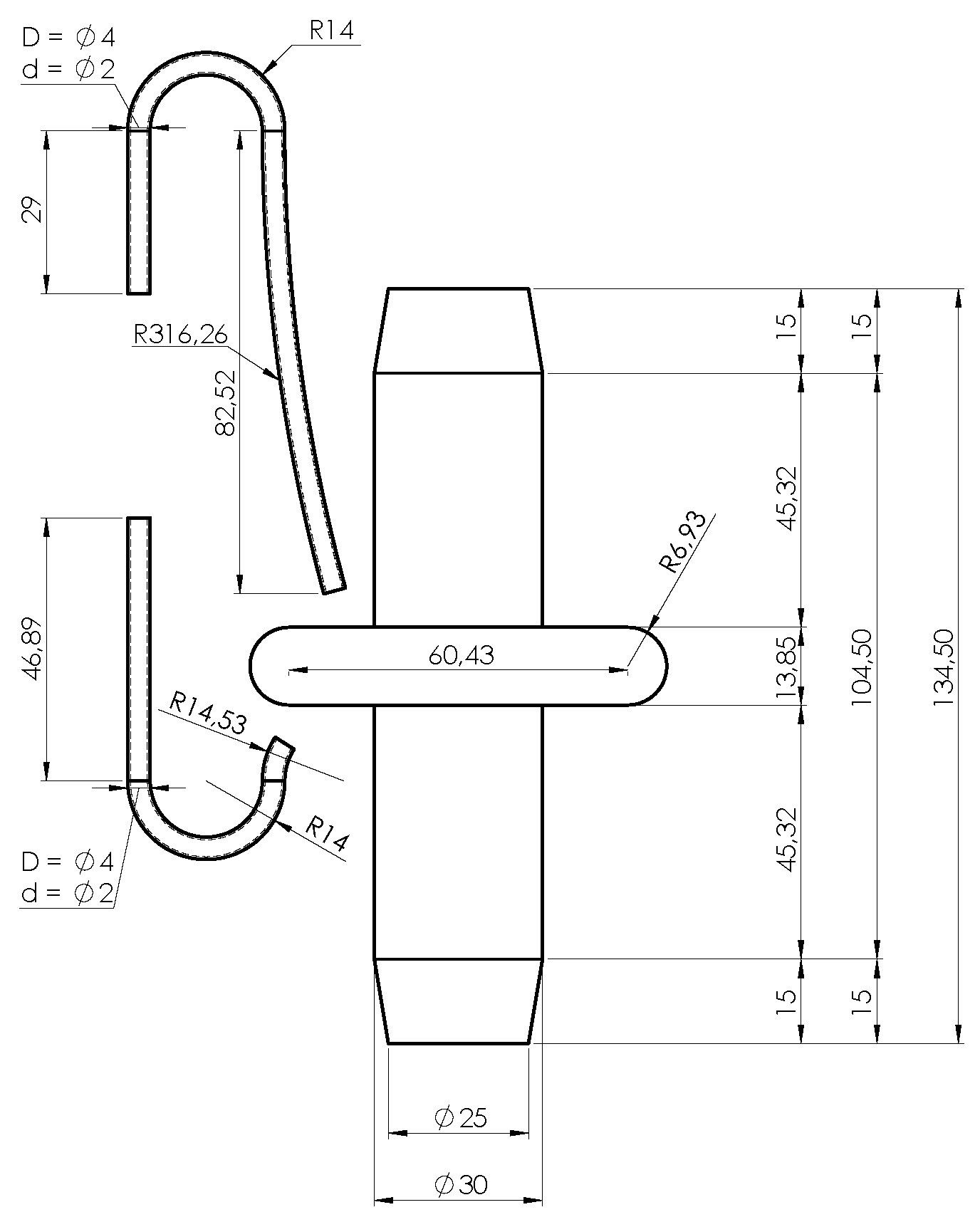}
 \caption{Sketch of the flanged shaft}
 \label{fig:flangesketch}
\end{figure}

As a second test case, we consider the cooling of a flanged shaft by
cold high pressured air, a process that's also known as gas
quenching. The complete process consists of the inductive heating of a
steel rod, the forming of the hot rod into a flanged shaft, a
transport to a cooling unit and the cooling process. Here, we consider
only the cooling, meaning that we have a hot flanged shaft that is
cooled by cold high pressured air coming out of small tubes. We
consider a two dimensional cut through the domain and assume symmetry
along the horizontal axis, resulting in one half of the flanged shaft
and two tubes blowing air at it, see figure \ref{fig:flangesketch}. 
We assume that the air leaves the tube in a straight and uniform way at
a Mach number of 1.2. Furthermore, we assume a freestream in $x$-direction of Mach 0.005. This is mainly to avoid numerical difficulties at Mach 0, but could model a draft in the workshop. The Reynolds number is  $Re=2500$ and the Prandtl number $Pr=0.72$. 
\begin{figure}[ht]
  \centering
  \subfigure[Entire mesh]{
    \includegraphics[width=0.5
    \linewidth]{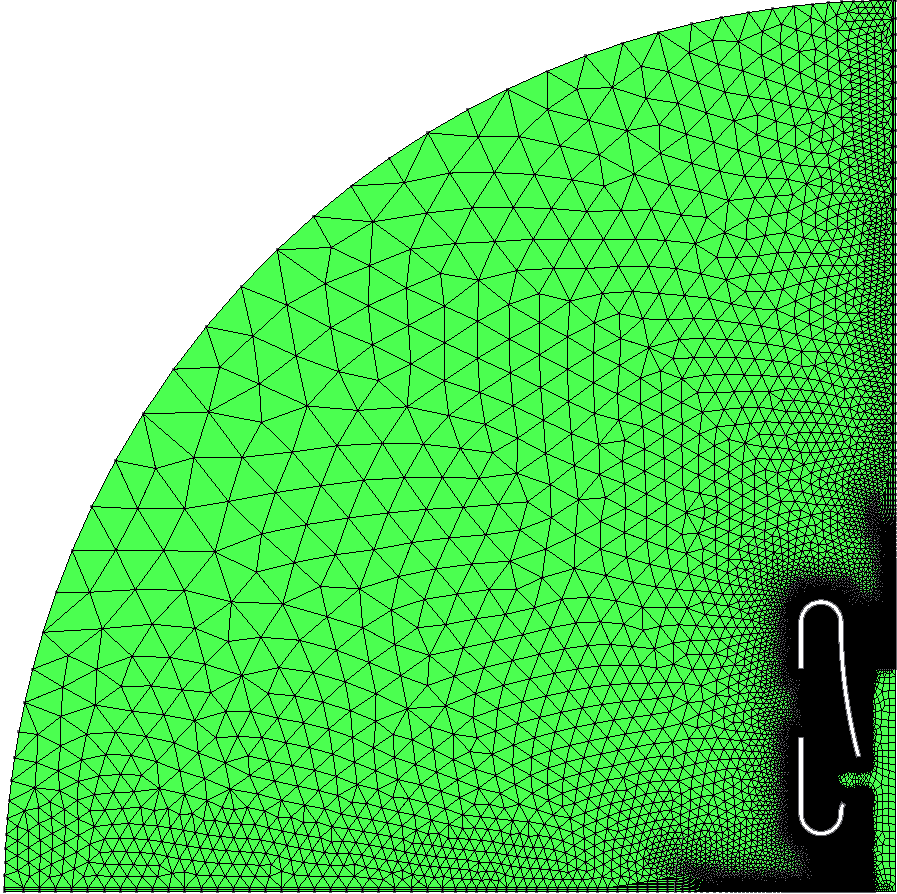}}
\hfill
  \subfigure[Mesh zoom]{
    \includegraphics[width=0.45 \linewidth]{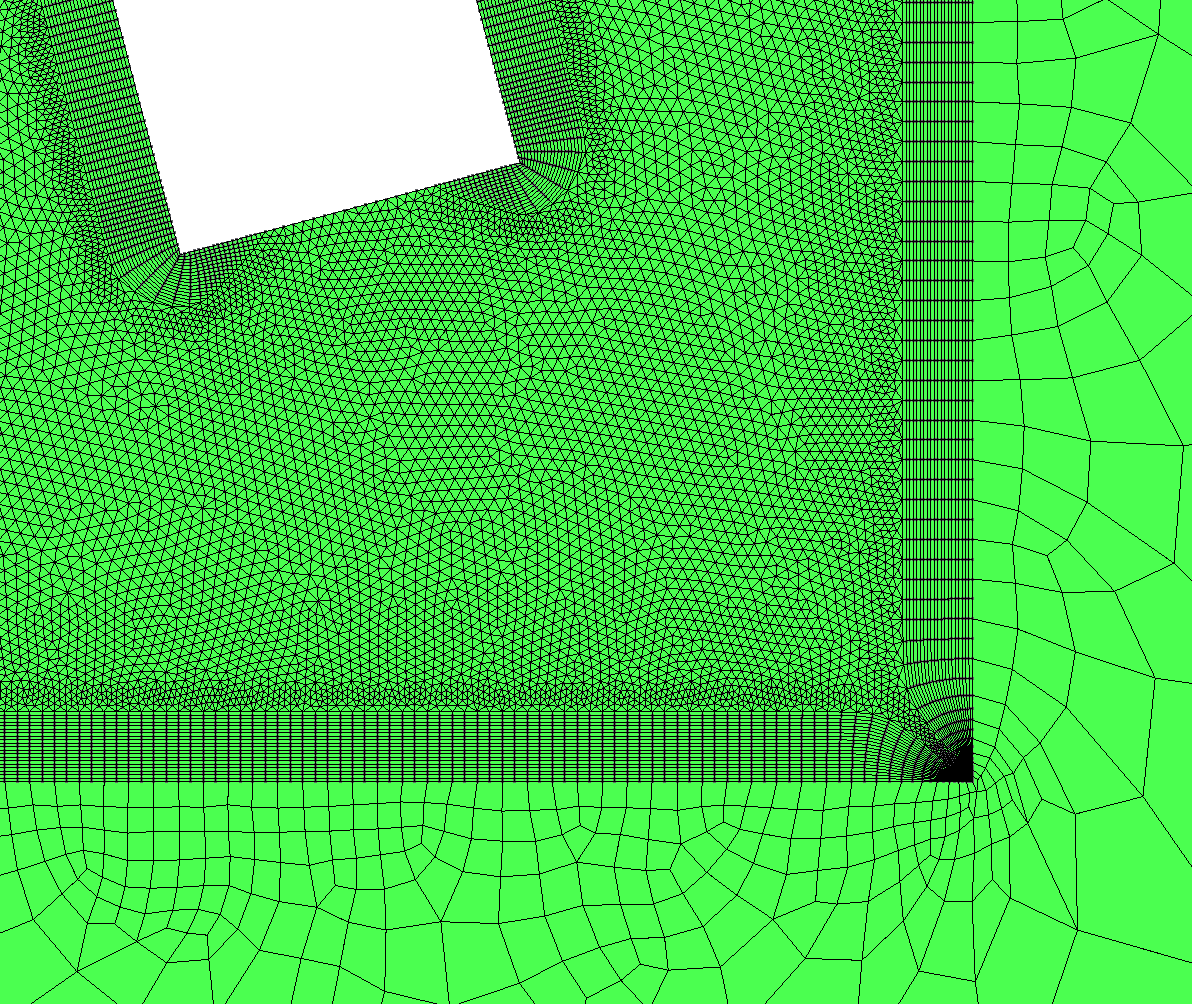}}
  \caption{Full grid (left) and zoom into shaft region (right)}
  \label{fig:flangegrid}
\end{figure}

The grid consists of 279212 cells in the fluid, which is the dual grid of an
unstructured grid of quadrilaterals in the boundary layer and
triangles in the rest of the domain, and 1997 quadrilateral elements in the structure. It is illustrated in figure \ref{fig:flangegrid}. 

\begin{figure}[h!]
  \centering
  \subfigure{
   \includegraphics[width=0.47 \linewidth]{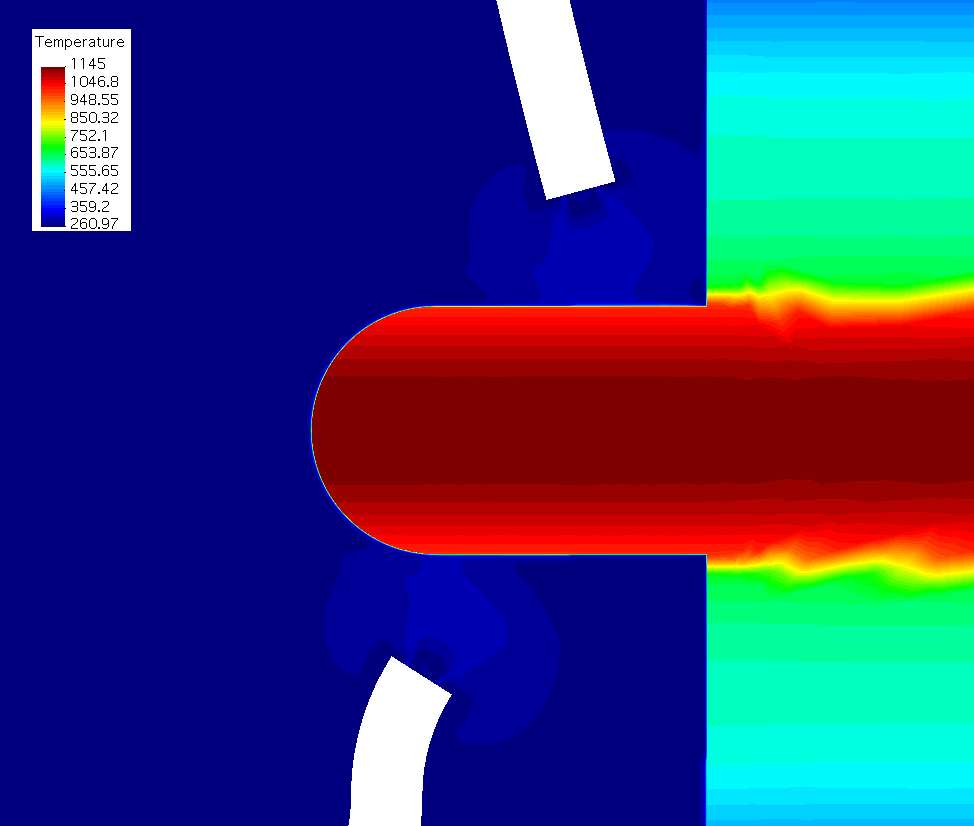}}
\hfill
\subfigure{
   \includegraphics[width=0.47 \linewidth]{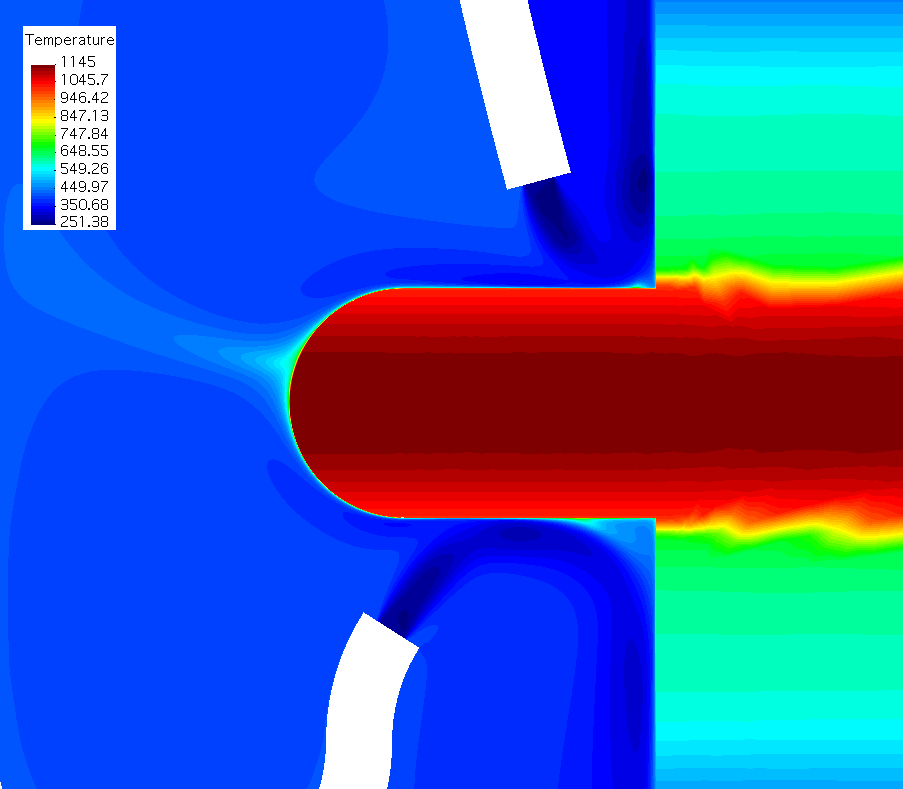}}
 \caption{Temperature distribution in fluid and structure at $t=0s$ (left) and $t=1s$ (right).}
  \label{fig:temperatureflange}
%
 \end{figure}

To obtain initial conditions for the subsequent tests, we use the
following procedure: We define a first set of initial conditions by
setting the flow velocity to zero throughout and choose the structure
temperatures at the boundary points to be equal to temperatures that
have been measured by a thermographic camera. Then, setting the $y$-axis on the symmetry axis of the flange, we set the temperatur at each
horizontal slice to the temperature at the correspoding boundary
point. Finally, to
determine the actual initial conditions, we
compute $10^{-5}$ seconds of real time using the coupling solver with
a fixed time step size of $\Delta t=10^{-6}s$. This means, that the
high pressured air is coming out of the tubes and the first front has
already hit the flanged shaft. This solution is illustrated in figure~\ref{fig:temperatureflange} (left). The wiggles in the structure are due to visualization artifacts. 

Now, we compute 1 second of real time using the time adaptive algorithm with different tolerances and an initial time step size of $\Delta t=10^{-6}s$. This small initial step size is necessary to prevent instabilities in the fluid solver. During the course of the computation, the time step size is increased until it is on the order of $\Delta t=0.1s$, which demonstrates the advantages of the time adaptive algorithm and reaffirms that it is this algorithm that we need to compare to. In total, the time adaptive method needs 22, 41, 130 and 850 time steps to reach $t=1s$ for the different tolerances, compared to the $10^6$ steps the fixed time step method would need. The solution at the final time is depicted in figure~\ref{fig:temperatureflange} (right). As can be seen, the stream of cold air is deflected by the shaft. 

We then compare the total number of iterations for the different vector extrapolation methods, see table \ref{tab:flangedvectorextr}. As before, the vector extrapolation methods have almost no effect on the number of iterations. 

\begin{table}[h!]
\caption{\label{tab:flangedvectorextr}Total number of iterations for 1 sec of real time for different vector extrapolation methods}
\centering
\begin{tabular}{*{6}{c}}
\hline
$TOL$ & & No relax. & Aitken & MPE & RRE \\
\hline
$10^{-2}$ &$\#$Iterations &52 &52 &52  &52\\
$10^{-3}$ &$\#$Iterations &127  &128  &127  &127 \\
$10^{-4}$ &$\#$Iterations &433   &430  &433  &433 \\
$10^{-5}$ &$\#$Iterations &2874   &2859  &2874  &2874 \\
\hline
\end{tabular}
\end{table}

Finally, we consider extrapolation based on the time integration scheme. In table \ref{tab:extrap-shaft}, the total number of iterations for 1 second of real time is shown. As before, the extrapolation methods cause a noticable decrease in the total number of fixed point iterations, with linear extrapolation performing better than the quadratic version. The speedup from linear extrapolation is between 20\% and 30\%, compared to the results obtained without extrapolation. 
\begin{table}[h!]
\caption{\label{tab:extrap-shaft}Total number of iterations for 1 sec of
  real time for different extrapolation methods in time}
\centering
\begin{tabular}{*{4}{c}}
\hline
$TOL$ & none & lin. & quad. \\
\hline
$10^{-2}$ &52 &42 &47 \\
$10^{-3}$ &127 &97 &99 \\
$10^{-4}$ &433 &309 &312 \\
$10^{-5}$ &2874 &1805 &1789 \\
\hline
\end{tabular}
\end{table}


\section{Summary and Conclusions}

We considered a time dependent thermal fluid structure interaction problem where a
nonlinear heat equation to model steel is coupled with the
compressible Navier-Stokes equations. The coupling is performed in a Dirichlet-Neumann manner. As a fast base solver, a higher order
time adaptive method is used for time integration. This method is significantly more efficient than a fixed time step method and is therefore the scheme to beat. 

To reduce the
number of fixed point iterations in a partitioned spirit, first different vector
extrapolation techniques, namely Aitken Relaxation, MPE and RRE were
compared. These have a negligible effect, since they are only useful
when a large number of iterations is needed per system and the time adaptive method results in only two iterations being necessary per time step. However, extrapolation based on the time integration
method works from the first iteration and reduces the number of
iterations by up to $50 \%$. Hereby, linear extrapolation works better than quadratic. 

The combined time adaptive method with linear extrapolation thus allows to solve real life problems at engineering tolerances using only a couple dozen calls to the fluid and structure solver.

\section*{Acknowledgement}

Financial support has been provided by the German Research Foundation (DFG) via the Sonderforschungsbereich Transregio 30, projects C1 and C2.


\begin{thebibliography}{10}

\bibitem{arnold:10}
{\sc M.~Arnold}, {\em {Stability of Sequential Modular Time Integration Methods
  for Coupled Multibody System Models}}, J. Comput. Nonlinear Dynam., 5 (2010),
  pp.~1--9.

\bibitem{banka:05}
{\sc A.~L. Banka}, {\em {Practical Applications of CFD in heat processing}},
  Heat Treating Progress,  (2005).

\bibitem{birken:14}
{\sc P.~Birken}, {\em {Termination criteria for inexact fixed point schemes}},
  Numer. Linear Algebra Appl., submitted.

\bibitem{birkenhabil}
{\sc P.~Birken}, {\em {Numerical Methods for the Unsteady Compressible
  Navier-Stokes Equations}}, Habilitation Thesis, University of Kassel, 2012.

\bibitem{biquhm:10}
{\sc P.~Birken, K.~J. Quint, S.~Hartmann, and A.~Meister}, {\em {Chosing norms
  in adaptive FSI calculations}}, PAMM, 10 (2010), pp.~555--556.

\bibitem{biquhm:11}
\leavevmode\vrule height 2pt depth -1.6pt width 23pt, {\em {A Time-Adaptive
  Fluid-Structure Interaction Method for Thermal Coupling}}, Comp. Vis. in
  Science, 13 (2011), pp.~331--340.

\bibitem{buchli:10}
{\sc J.~M. Buchlin}, {\em {Convective Heat Transfer and Infrared
  Thermography}}, J. Appl. Fluid Mech., 3 (2010), pp.~55--62.

\bibitem{erbdue:12}
{\sc P.~Erbts and A.~D\"{u}ster}, {\em {Accelerated staggered coupling schemes
  for problems of thermoelasticity at finite strains}}, Comp. \& Math. with
  Appl., 64 (2012), pp.~2408--2430.

\bibitem{farhat:04}
{\sc C.~Farhat}, {\em {CFD-based Nonlinear Computational Aeroelasticity}}, in
  Encyclopedia of Computational Mechanics, E.~Stein, R.~de~Borst, and T.~J.~R.
  Hughes, eds., vol.~3: Fluids, John Wiley \& Sons, 2004, ch.~13, pp.~459--480.

\bibitem{gerfeg:97}
{\sc T.~Gerhold, O.~Friedrich, J.~Evans, and M.~Galle}, {\em {Calculation of
  Complex Three-Dimensional Configurations Employing the DLR-TAU-Code}}, AIAA
  Paper, 97-0167 (1997).

\bibitem{giles:97}
{\sc M.~B. Giles}, {\em {Stability Analysis of Numerical Interface Conditions
  in Fluid-Structure Thermal Analysis}}, Int. J. Num. Meth. in Fluids, 25
  (1997), pp.~421--436.

\bibitem{hefiba:01}
{\sc U.~Heck, U.~Fritsching, and K.~Bauckhage}, {\em {Fluid flow and heat
  transfer in gas jet quenching of a cylinder}}, International Journal of
  Numerical Methods for Heat \& Fluid Flow, 11 (2001), pp.~36--49.

\bibitem{hinrad:06}
{\sc M.~Hinderks and R.~Radespiel}, {\em {Investigation of Hypersonic Gap Flow
  of a Reentry Nosecap with Consideration of Fluid Structure Interaction}},
  AIAA Paper, 06-1111 (2006).

\bibitem{kuewal:08}
{\sc U.~K\"{u}ttler and W.~A. Wall}, {\em {Fixed-point fluid–structure
  interaction solvers with dynamic relaxation}}, Comput. Mech., 43 (2008),
  pp.~61--72.

\bibitem{letmou:01}
{\sc P.~{Le Tallec} and J.~Mouro}, {\em {Fluid structure interaction with large
  structural displacements}}, Comp. Meth. Appl. Mech. Engrg., 190 (2001),
  pp.~3039--3067.

\bibitem{lior:04}
{\sc N.~Lior}, {\em {The cooling process in gas quenching}}, J. Materials
  Processing Technology, 155-156 (2004), pp.~1881--1888.

\bibitem{mehta:05}
{\sc R.~C. Mehta}, {\em {Numerical Computation of Heat Transfer on Reentry
  Capsules at Mach 5}}, AIAA-Paper 2005-178,  (2005).

\bibitem{mivabo:06}
{\sc C.~Michler, E.~H. van Brummelen, and R.~de~Borst}, {\em
  {Error-amplification Analysis of Subiteration-Preconditioned GMRES for
  Fluid-Structure Interaction}}, Comp. Meth. Appl. Mech. Eng., 195 (2006),
  pp.~2124--2148.

\bibitem{olssoe:98}
{\sc H.~Olsson and G.~S\"{o}derlind}, {\em {Stage value predictors and
  efficient Newton iterations in implicit Runge-Kutta methods}}, SIAM J. Sci.
  Comput., 20 (1998), pp.~185--202.

\bibitem{quhrss:11}
{\sc K.~J. Quint, S.~Hartmann, S.~Rothe, N.~Saba, and K.~Steinhoff}, {\em
  {Experimental validation of high-order time integration for non-linear heat
  transfer problems}}, Comput. Mech., 48 (2011), pp.~81--96.

\bibitem{sidi:12}
{\sc A.~Sidi}, {\em {Review of two vector extrapolation methods of polynomial
  type with applications to large-scale problems}}, J. Comp. Phys., 3 (2012),
  pp.~92--101.

\bibitem{stshle:06}
{\sc P.~Stratton, I.~Shedletsky, and M.~Lee}, {\em {Gas Quenching with
  Helium}}, Solid State Phenomena, 118 (2006), pp.~221--226.

\bibitem{zubobi:07}
{\sc S.~van Zuijlen, A. H.~Bosscher and H.~Bijl}, {\em {Two level algorithms
  for partitioned fluid–structure interaction computations}}, Comp. methods
  appl. mech. eng.,  (2007).

\bibitem{viladv:06}
{\sc J.~Vierendeels, L.~Lanoye, J.~Degroote, and P.~Verdonck}, {\em {Implicit
  Coupling of Partitioned Fluid-Structure Interaction Problems with Reduced
  Order Models}}, Comp. \& Struct., 85 (2007), pp.~970--976.

\bibitem{wesast:07}
{\sc U.~Weidig, N.~Saba, and K.~Steinhoff}, {\em {Massivumformprodukte mit
  funktional gradierten Eigenschaften durch eine differenzielle
  thermo-mechanische Prozessf\"{u}hrung}}, WT-Online,  (2007), pp.~745--752.

\bibitem{zienkiewicztaylor00a}
{\sc O.~Zienkiewicz and R.~Taylor}, {\em The {F}inite {E}lement {M}ethod},
  Butterworth Heinemann, 2000.

\end{thebibliography}
\end{document}